\documentclass[reqno]{amsart}
\usepackage{a4wide,calrsfs}
\usepackage{amsmath,bbm}
\usepackage{color}
\usepackage{epsfig}
\usepackage{subfigure}
\usepackage{amssymb}
\usepackage{enumerate}
\usepackage{graphicx}
\usepackage{subfigure}
\usepackage{amsmath}
\usepackage{rotating}
\usepackage{stmaryrd}
\usepackage{upgreek}

\newtheorem{theorem}{Theorem}[section]

\allowdisplaybreaks
\numberwithin{equation}{section}

\title[A proof of the conjectures by Yu-Dong Wu and H.M. Srivastava]{
A proof of the three geometric inequalities conjectured by Yu-Dong Wu and H.M. Srivastava}

\author[Coronel]{An{\'\i}bal\ Coronel$^\dag$}
\author[Huancas]{Fernando Huancas$^\ddag$}

\thanks{$^\dag$ GMA, Departamento de Ciencias B\'asicas,
Facultad de Ciencias, Universidad del B\'{\i}o-B\'{\i}o,
Campus Fernando May, Chill\'{a}n, Chile,
E-mail: {\tt acoronel@ubiobio.cl}}

\thanks{$^\ddag$ GMA and Doctorado en Matem\'atica Aplicada,
Facultad de Ciencias, Universidad del B\'{\i}o-B\'{\i}o,
Campus Fernando May, Chill\'{a}n, Chile,
E-mail: {\tt fihuanca@gmail.com}}

\date{\today}

\begin{document}

\begin{abstract}
In this short note the authors give answers to the three open problems 
formulated by
Wu and Srivastava [{\it Appl. Math. Lett. 25 (2012), 1347--1353}].
We disprove the Problem 1,
by showing that there exists a triangle
which does not satisfies the proposed inequality.
We prove the inequalities conjectured in Problems 2 and 3. 
Furthermore, we introduce an optimal refinement of the
inequality conjectured on Problem 3.
\end{abstract}


\maketitle

\section{Introduction and conjectures of Yu-Dong Wu and H.M. Srivastava }
The geometric inequalities are relevant in several areas of the science and engineering
\cite{cvetkovski2012inequalities,cloud1998inequalities,lieb2002inequalities,
MR1022443,MR1785857,Srivastava20112349,Wu2010761}. 
The methodologies to prove the geometric inequalities is disperse
see for instance \cite{satnoianu2001,pecaric_1997,yanglu1999,yudong2013}.
In a broad sense, there exist some of methodologies which are
based on analytical methods, other in integral and differential calculus,
and other on geometric methods.  The methodology
of this paper, in spite of the basic arguments,
can be considered belongs to the analytical methods since 
the results are strongly dependent 
on the  analytical methodology introduced in~\cite{MR1785857,wusrivastava_2012}.

The focus of this short note is the open problems given in~\cite{MR1785857,wusrivastava_2012}.
Indeed, we introduce some notation and then we recall the conjectures.
Let us consider a triangle $\triangle ABC$ with angles $A,B$ and $C$, we denote
by $a,b,c,s$ and $r$, the lengths of the corresponding opposite sides, the semiperimeter
and the inradius, respectively. Then, using the symbol $\sum$ to denote a  cyclic sum, i. e.
\begin{eqnarray*}
\sum f(b,c)=f(a,b)+f(b,c)+f(c,a),
\end{eqnarray*}
we have that, the following geometric inequality
\begin{eqnarray}
&&2\sqrt{2}\; s
\le 
	\sum \sqrt{a^2+b^2}
<
	(2+\sqrt{2})\; s-3\sqrt{3}(2-\sqrt{2})r,
\label{eq:wu_srivastava}
\end{eqnarray}
holds. In (\ref{eq:wu_srivastava}), the left inequality  can be proved by the
Power-Mean Inequality and the right inequality was 
recently proved by Wu and Srivastava in \cite{wusrivastava_2012}.
It was originally proposed
by Wu \cite{wu_2001}, inspired in the following inequality \cite{wu_2001,zhou_1996}:
\begin{eqnarray}
\sum \sqrt{a^2+b^2}<(2+\sqrt{2})s.
\label{eq:zhouhu}
\end{eqnarray}
Moreover, it is known the following generalized version of
(\ref{eq:zhouhu}):  
\begin{eqnarray}
\sum \sqrt[n]{a^n+b^n}
<
(2+\sqrt[n]{2})s,
\qquad n\in\mathbb{N}-\{1\},
\label{eq:ye}
\end{eqnarray}
holds true for $\mathbb{N}$ the set of positive integers,
see \cite{ye_2003}. 
Then, by analogy to the extension of Ye \cite{ye_2003},
Wu and Srivastava propposed a generalization of  the 
right inequality in \cite{wusrivastava_2012}. More precisely, they defined
the following conjecture:
\begin{enumerate}
\item[]{\bf Conjecture 1.} {\it For a given triangle $\triangle ABC$, 
if $n\in\mathbb{N}-\{1,2\}$, then prove
or disprove the following inequality:}
\begin{eqnarray}
	\sum \sqrt[n]{a^n+b^n}
<
	(2+\sqrt[n]{2})s-3\sqrt{3}(2-\sqrt[n]{2})r.
\label{eq:conjecture_wu_srivastava}
\end{eqnarray}
\end{enumerate}
Aditionally, in \cite{wusrivastava_2012}
Wu and Srivastava conjecture the following two inequalities: 
\begin{enumerate}
\item[]{\bf Conjecture 2.} {\it Let $a_i$ ($i=1,\ldots,6$) denote the lengths
of the edges of a given tetrahedron $ABCD$. Also let $\rho$ be the inradius
of the  tetrahedron. Then, determine the best constant $k$ such the following
inequality holds true:}
\begin{eqnarray}
	\sum_{1\le i,j\le 6} \sqrt{a^2_i+b^2_j}
\le k\sum_{i=1}^6 a_i.
\label{eq:conjecture_wu_srivastava:2}
\end{eqnarray}

\item[]{\bf Conjecture 3.} {\it 
Let us denote by $k_0$ denotes the best constant $k$ for the inequality 
(\ref{eq:conjecture_wu_srivastava:2})
for a given tetrahedron $ABCD$. Then, prove
or disprove the following inequality }
\begin{eqnarray}
	\sum_{1\le i,j\le 6} \sqrt{a^2_i+b^2_j}
\le k_0\sum_{i=1}^6 a_i-6\sqrt{6}(2k_0-5\sqrt{2})\rho.
\label{eq:conjecture_wu_srivastava:3}
\end{eqnarray}
\end{enumerate}
We note that,
the inequality (\ref{eq:conjecture_wu_srivastava}) looks as a nice 
generalization of (\ref{eq:wu_srivastava}). However, unfortunately,
we show that the inequality (\ref{eq:conjecture_wu_srivastava})
is not always true, see subsection~\ref{sec:conj1}. 
In subsection ~\ref{sec:conj2}, we determine that
(\ref{eq:conjecture_wu_srivastava:2}) holds true with 
$k=2+\sqrt{2}.$
Furthermore, denoting by $k_0=2+\sqrt{2}$, in subsection ~\ref{sec:conj3},
we prove that
the following inequality 
\begin{eqnarray}
	\sum_{1\le i,j\le 6} \sqrt{a^2_i+b^2_j}
\le k_0\sum_{i=1}^6 a_i-12\sqrt{3}(2-\sqrt{2})\rho,
\label{eq:conjecture_wu_srivastava:3_true}
\end{eqnarray}
which is a refinement of (\ref{eq:conjecture_wu_srivastava:3}).
Then, the Conjecture 3 holds true but the inequality is not optimal.
Hence, sumarizing the contribution of this
sort note we have the following theorem:
\begin{theorem}
 \label{teo:main}
 Let $a_i$ ($i=1,\ldots,6$) denote the lengths
of the edges of a given tetrahedron $ABCD$. Also let $\rho$ be the inradius
of the  tetrahedron. Then, the inequalities 
(\ref{eq:conjecture_wu_srivastava:2}) and
(\ref{eq:conjecture_wu_srivastava:3_true}) holds true with $k=k_0=2+\sqrt{2}.$
\end{theorem}

\section{Proofs of Conjectures}

\subsection{Conterexample for Conjecture 1}
\label{sec:conj1}
We given a Conterexample which proves that (\ref{eq:conjecture_wu_srivastava}) 
does not holds true. Indeed, let us consider 
$n=3$ and the right triangle $\triangle ABC$ with
$a=3,\;b=1 $ and $c=\sqrt{10}$.
Then, $s=(4+\sqrt{10})/2$, $r=3/(4+\sqrt{10})$ and
clearly the inequality (\ref{eq:conjecture_wu_srivastava}) is reversed, since
\begin{eqnarray*}
\sum \sqrt[3]{a^3+b^3}&\approx& 10.116536541585731
\quad\mbox{and}\quad
\\
(2+\sqrt[3]{2})s-3\sqrt{3}(2-\sqrt[3]{2})r&\approx& 10.063472825231253.
\end{eqnarray*}
Furthermore,  if for the given triangle we define the 
function $g:\mathbb{N}-\{1\}\to\mathbb{R}$ as follows
\begin{eqnarray*}
g(n)=(2+\sqrt[n]{2})s-3\sqrt{3}(2-\sqrt[n]{2})r-\sum \sqrt[n]{a^n+b^n},
\end{eqnarray*}
we note that $g$ is strictly decreasing with $g(2)>0>g(3)$. 
Thus, for the given right triangle the inequality
(\ref{eq:conjecture_wu_srivastava}) is false for all $n\in \mathbb{N}-\{1,2\}$.

\subsection{Proof of Conjecture 2}
\label{sec:conj2}

We apply the inequality (\ref{eq:zhouhu}) 
over each face of the tetrahedron $ABCD$
and, naturally, we get the following optimal 
estimates
\begin{eqnarray*}
&&\sum_{1\le i,j\le 6} \sqrt{a^2_i+b^2_j}
\\
&&
\qquad=
\Bigg[\sqrt{a^2_1+a^2_3}+\sqrt{a^2_3+a^2_2}+\sqrt{a^2_2+a^2_1}\Bigg]
+\Bigg[\sqrt{a^2_3+a^2_4}+\sqrt{a^2_4+a^2_5}+\sqrt{a^2_5+a^2_3}\Bigg]
\\
&&
\qquad\;\;
+\Bigg[\sqrt{a^2_1+a^2_4}+\sqrt{a^2_4+a^2_6}+\sqrt{a^2_6+a^2_1}\Bigg]
+\Bigg[\sqrt{a^2_6+a^2_5}+\sqrt{a^2_5+a^2_2}+\sqrt{a^2_2+a^2_4}\Bigg]
\\
&&
\qquad
\le
\left(1+\frac{\sqrt{2}}{2}\right)(a_1+a_2+a_3)
+\left(1+\frac{\sqrt{2}}{2}\right)(a_3+a_4+a_5)
\\
&&
\qquad\;\;
+\left(1+\frac{\sqrt{2}}{2}\right)(a_1+a_4+a_6)
+\left(1+\frac{\sqrt{2}}{2}\right)(a_2+a_5+a_6)
\\
&&
\qquad
\le
(2+\sqrt{2})\sum_{i=1}^6 a_i.
\end{eqnarray*}
Then, we have that the inequality (\ref{eq:conjecture_wu_srivastava:2}) holds with optimal
constant $k=2+\sqrt{2}$.

\subsection{Proof of Conjecture 3}
\label{sec:conj3}

Let us denote by $k_0=2+\sqrt{2}$ 
and by $r_i$ for $i=1,\ldots,4$ the inradius of the faces of the
tetrahedron $ABCD$. Then,  
applying the inequality (\ref{eq:conjecture_wu_srivastava})
over each face of the tetrahedron $ABCD$
we fid that (\ref{eq:conjecture_wu_srivastava:3_true}) holds
true, since
\begin{eqnarray*}
\sum_{1\le i,j\le 6} \sqrt{a^2_i+b^2_j}
&\le& k_0\sum_{i=1}^6 a_i-3\sqrt{3}(2-\sqrt{2})\sum_{i=1}^4 r_i
\\
&\le& k_0\sum_{i=1}^6 a_i-12\sqrt{3}(2-\sqrt{2})\rho.
\end{eqnarray*}
Moreover, we note that
\begin{eqnarray*}
 k_0\sum_{i=1}^6 a_i-12\sqrt{3}(2-\sqrt{2})\rho
 \le k_0\sum_{i=1}^6 a_i-6\sqrt{6}(2k_0-5\sqrt{2})\rho.
\end{eqnarray*}
Then, the inequality (\ref{eq:conjecture_wu_srivastava:3}) is true, 
but it is not optimal.

\section*{Acknowledgement}
We acknowledge the support of ``Univesidad del B{\'\i}o-B{\'\i}o" (Chile) through
the research projects 124109 3/R, 104709 01 F/E and 
121909 GI/C.

\end{document}